\documentclass[12pt]{article}

\usepackage[a4paper]{geometry}
\usepackage{graphicx}  % Add graphics capabilities
\usepackage{amsmath}  % Better maths support & more symbols
\usepackage{amssymb,amsthm}
\usepackage{color}
\usepackage{tikz}
\usepackage{pgfplots}
\usepackage{hyperref}

\newcommand{\R}{{\mathbb R}}\newcommand{\N}{{\mathbb N}}
\newcommand{\T}{{\mathbb T}}
\newcommand{\Z}{{\mathbb Z}}

\let\epsilon\varepsilon
\let\theta\vartheta
\let\hat\widehat

\newtheorem{theorem}{Theorem}[section]\newtheorem{lemma}[theorem]{Lemma}
\newtheorem{definition}[theorem]{Definition}

\newtheorem{corollary}[theorem]{Corollary}

\newtheorem{remark}[theorem]{Remark}

\title{On the KdV approximation for a Boussinesq equation posed  on
the infinite necklace graph}
\author{Wolf-Patrick D\"ull, Guido Schneider, Raphael Taraca \\ 
{\small  Institut f\"{u}r Analysis, Dynamik und Modellierung,
Universit\"{a}t Stuttgart,} \\
{\small Pfaffenwaldring 57, D-70569 Stuttgart, Germany }  } 

\begin{document}

\maketitle

\begin{abstract}
We consider a Boussinesq equation posed
on the infinite periodic necklace graph.
For the description of long wave traveling waves we derive the KdV equation 
and establish the validity of this formal approximation by providing estimates for the error. 
The proof is based on  suitable energy estimates.
\end{abstract}

%{\color{red}     
%}

\section{Introduction}

The KdV equation  
\begin{equation} \label{eqkdv}
\partial _{T}A=\nu_1 \partial ^{3}_{X}A+\nu_2 \partial _{X}(A^{2}),
\end{equation}
with $ T ,X,A(X,T) ,\nu_1,\nu_2 \in \R $ is well known 
for being a completely integrable Hamiltonian system possessing solitons, cf. \cite{Drazin}. 
It can be derived as an approximation equation through a multiple scaling perturbation ansatz for various original systems. In this way, 
it 
appears as a universal amplitude equation 
describing slow modulations in time and space of long wave solutions
for many dispersive systems. 

For spatially homogeneous systems
the last decades saw various approximation results guaranteeing 
that  such formally derived KdV equations make correct predictions about the dynamics of the 
original systems. 
Such error estimates have been shown for instance in \cite{Cr85,KN86,SW00,SW02,Du12} 
for the description of surface water waves or in \cite{SW00equa}
for the  FPUT model.

For spatially periodic systems which  are not perturbations of a spatially homogeneous system
only a few results exist. In \cite{CCPS12,GMWZ14} the validity of the KdV approximation for poly-atomic FPUT problems
have been established. In \cite{BDS19} such estimates have been proven for a Boussinesq equation with spatially periodic  coefficients. 

In the last years quantum graphs attracted a lot of attention, cf. Section \ref{sec2}.
An introduction into the mathematical theory of quantum graphs can be found for instance in \cite{Kuchment}. They are used as models for nanotechnological devices such as nanotubes or 
graphene. For the understanding of the qualitative behavior of solutions of the 
original systems posed on infinitely extended periodic quantum graphs
the derivation of approximation equations 
through multiple scaling perturbation theory can be used again.
In \cite{GPS16,GSU20} the NLS approximation and the Dirac approximation 
have been justified  for 1D and 2D infinite periodic quantum graphs. 
In this paper 
we are interested in the question whether the KdV equation can be justified 
as a modulation equation for the dynamics of  the Boussinesq equation posed on 
the infinite periodic  necklace graph. 

Before we explain subsequently what is meant by
posing the Boussinesq equation on the infinite periodic necklace  graph $\Gamma$ 
we first
consider the Boussinesq equation 
\begin{equation}\label{Boussintro}
\partial_t^2 u = \partial_x^2 u + \partial_t^2 \partial_x^2 u + \partial_x^2 (u^2),
\end{equation}
on the real line, i.e. $ x \in \R $, with $ t ,u(x,t) \in \R $
and explain in the following remark how the KdV approximation is obtained in this situation.
\begin{remark}{\rm
For $ x \in \R $ by inserting the ansatz 
\begin{equation}
\label{kdvhomoansatz}
 u(x,t) = \varepsilon^2 A(X,T) ,
 \end{equation}
with 
$
X = \varepsilon (x-c t)  $ and $ T = \varepsilon^3 t $
into \eqref{Boussintro} we obtain 
$$ 
\varepsilon^4 c^2 \partial_X^2 A  - 2 \varepsilon^6 c \partial_X \partial_T A + \mathcal{O}(\varepsilon^8) 
= \varepsilon^4 \partial_X^2 A +  c^2 \varepsilon^6  \partial_X^4 A + \mathcal{O}(\varepsilon^8) + \varepsilon^6 \partial_X^2 (A^2).
$$ 
Equating the coefficient of $ \varepsilon^4 $ to zero  gives $ c^2 = 1 $, and  
equating the coefficient of $ \varepsilon^6 $ to zero  gives 
the KdV equation 
\begin{equation}\label{kdvintro}
- 2 c  \partial_T A = c^2   \partial_X^3 A + \partial_X (A^2).
\end{equation}
The following approximation result holds, cf. \cite[\S 12]{SU17book}.

{\bf Theorem.}
Let $ A \in C([0,T_0],H^4) $ be a solution of the KdV equation \eqref{kdvintro}.
Then there exist $ C > 0 $ and $ \varepsilon_0 > 0 $ such that for all 
$ \varepsilon \in (0,\varepsilon_0) $ there are solutions $ u $ of the Boussinesq equation 
\eqref{Boussintro} with 
$$
\sup_{t \in [0,T_0/\varepsilon^3]} \sup_{x \in \R} |u(x,t)-\varepsilon^2 A( \varepsilon (x-c t) ,\varepsilon^3 t) | \leq C \varepsilon^3.
$$}
\end{remark}
It is the main goal of this paper to transfer this result to the Boussinesq equation 
posed  on the infinite periodic necklace  graph, cf. Theorem \ref{Main theorem}. 
One of the major difficulties in doing so comes  from the  non-smoothness of the solutions at the vertices 
of the graph.

For this specific example we develop an analysis which can be used for the 
justification of the KdV approximation on one-dimensional infinitely extended periodic quantum graphs. 
In Section \ref{sec9} we present a few other quantum graphs for which the 
KdV approximation can be justified. There, we also discuss why a general justification
result for the KdV approximation on quantum graphs would be only of minor use.

The plan of the paper is as follows. In Section \ref{sec2} and Section \ref{sec3} we start by explaining what is meant by posing 
the Boussinesq equation  on the infinite periodic necklace  graph and 
prove  the local existence and uniqueness of solutions.
After that we recall the spectral properties of the linearized problem in Section \ref{sec4}  and 
derive 
the KdV equation in Section \ref{sec5} and Section \ref{sec6}.
We derive the equations for the error, estimate the residual terms and the error with the help of a suitable chosen 
energy  in Sections \ref{secresidual}-\ref{sec8}.
The paper is closed with a discussion about possible extensions of the result in Section \ref{sec9}.
\medskip

{\bf Notation.}
The Sobolev space $H^s(\R)$ of $s$-times weakly differentiable $ L^2 $-functions is 
equipped with the norm 
\[ \| u \|_{H^s} = \left( \sum_{k=0}^s \| \partial_x^k u \|_{L^2} \right)^\frac{1}{2} \]
or  equivalently with the norm
\[ \| u \|_{H^s} = \| (1+(\cdot)^2)^{s/2} \hat{u}(\cdot) \|_{L^2}, \]
if $ s \geq 0 $ is a real number,
where $\hat{u}(k) = (\mathcal{F}u)(k)$ denotes the Fourier transform of $u$. Possibly different constants are denoted by $C$, if they do not depend on the small perturbation parameter $0<\varepsilon\ll 1$.
\medskip

{\bf Acknowledgement.}
The paper is partially supported by the Deutsche For\-schungsgemeinschaft DFG through the SFB 1173 ''Wave phenomena''
under the grant 258734477.

\section{The necklace graph as a quantum graph}

\label{sec2}

\begin{figure}[htbp]%  figure placement: here, top, bottom, or page
\centering
  \setlength{\unitlength}{1.1cm}
\begin{picture}(12,6)(2,0)
\put(2,5){\line(1,0){1}}
\put(3.5,5){\circle{1}}
\put(4,5){\line(1,0){1}}
\put(5.5,5){\circle{1}}
\put(6,5){\line(1,0){1}}
\put(7.5,5){\circle{1}}
\put(8,5){\line(1,0){1}}
\put(9.5,5){\circle{1}}
\put(10,5){\line(1,0){1}}
\put(11.5,5){\circle{1}}
\put(12,5){\line(1,0){1}}
\put(2,2){\line(1,0){1}}
\put(3,2.5){\line(1,0){1}}
\put(3,1.5){\line(1,0){1}}
\put(4,2){\line(1,0){1}}
\put(5,2.5){\line(1,0){1}}
\put(5,1.5){\line(1,0){1}}
\put(6,2){\line(1,0){1}}
\put(7,2.5){\line(1,0){1}}
\put(7,1.5){\line(1,0){1}}
\put(8,2){\line(1,0){1}}
\put(9,2.5){\line(1,0){1}}
\put(9,1.5){\line(1,0){1}}
\put(10,2){\line(1,0){1}}
\put(11,2.5){\line(1,0){1}}
\put(11,1.5){\line(1,0){1}}
\put(12,2){\line(1,0){1}}
\put(6.2,5.3){$\Gamma_{n,0}$}
\put(7.2,5.8){$\Gamma_{n,+}$}
\put(7.2,4){$\Gamma_{n,-}$}
\put(6.2,2.3){$I_{n,0}$}
\put(7.2,2.8){$I_{n,+}$}
\put(7.2,1){$I_{n,-}$}
  \end{picture}
\caption{The infinite periodic necklace graph $\Gamma$ shown in the upper panel
 is of the form 
$
\Gamma = \oplus_{n \in \mathbb{Z}} \Gamma_n $, 
with 
$ \Gamma_n = \Gamma_{n,0} \oplus \Gamma_{n, + } \oplus
\Gamma_{n, -}$,
where the $\Gamma_{n,0}$ are the  horizontal links 
between the circles and the $\Gamma_{n, \pm}$  the upper and
lower semicircles, all of the same length $\pi$, for $n \in
\mathbb{Z}$. The part  $\Gamma_{n,0}$ is identified
isometrically with the interval $I_{n,0} = \left[ 2\pi n, 2\pi n +
  \pi \right]$ and the  $\Gamma_{n, \pm}$  with the intervals
$I_{n, \pm } = \left[ 2\pi n + \pi, 2\pi (n + 1) \right]$. See the lower
panel. For a function $u : \Gamma \to \mathbb{C}$, we denote the part
on the interval $I_{n,0}$  with $u_{n,0}$
and the parts on the intervals $I_{n, \pm}$  with $u_{n, \pm}$.  \label{fig2}
}
\end{figure}
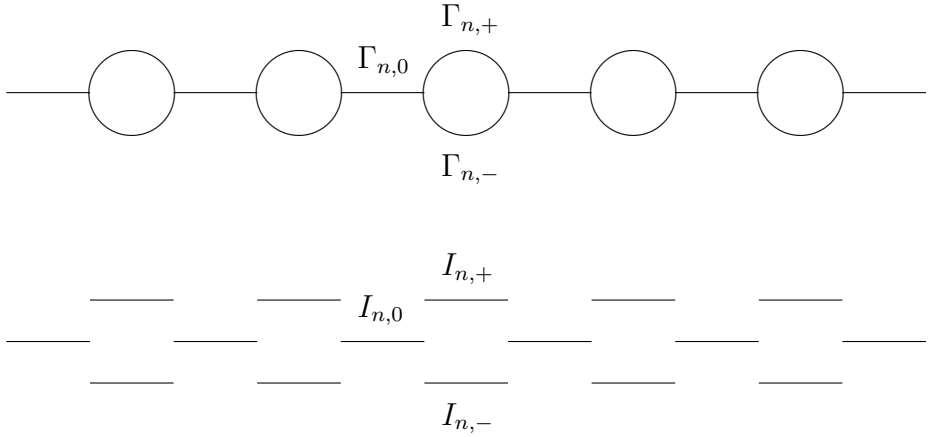

A quantum graph is a metric graph, i.e., a system of vertices and  connecting edges, equipped with the Laplace operator and corresponding boundary conditions at the vertices. The Laplace operator acts on functions defined along the edges.
In this paper we restrict ourselves to the infinite periodic necklace graph which  is sketched in Figure \ref{fig2}. 
In the caption of Figure  \ref{fig2} we  identify  this metric graph with a set of intervals connected at the vertices.

%\begin{definition}
On the edges of the necklace graph we consider the differential operator $ \partial_x^2 $ 
and at the vertices 
we assume so called Kirchhoff boundary conditions which are given by  the continuity of the functions, namely 
\begin{equation}\label{eq. KH1}
\begin{cases}
u_{n,0}(2\pi n+\pi,t)=u_{n,+}(2\pi n+\pi,t)=u_{n,-}(2\pi n+\pi,t), \\
u_{n+1,0}(2\pi(n+1),t)=u_{n,+}(2\pi(n+1),t)=u_{n,-}(2\pi(n+1),t)
\end{cases}
\end{equation}
and  the continuity of the fluxes 
\begin{equation}\label{eq. KH2}
\begin{cases}
\partial_x u_{n,0}(2\pi n+\pi,t)=\partial_x u_{n,+}(2\pi n+\pi,t)+\partial_x u_{n,-}(2\pi n+\pi,t), \\
\partial_x u_{n+1,0}(2\pi(n+1),t)=\partial_x u_{n,+}(2\pi(n+1),t)+\partial_x u_{n,-}(2\pi(n+1),t).
\end{cases}
\end{equation}
We collect the functions $ u_{n,0} $, $u_{n,-} $, and $u_{n,+} $ in 
$$
u_0(x, t )=
\begin{cases}
u_{n,0}(x, t ), & x\in I_{n,0} \\
0, & x\in I_{n,\pm}
\end{cases},
\quad n\in\Z
$$
and 
$$
u_\pm(x, t )=
\begin{cases}
u_{n,\pm}(x, t ), & x\in I_{n,\pm} \\
0, & x\in I_{n,0}
\end{cases},
\quad n\in\Z,
$$
with their supports  
\[ I_0 = \bigcup_{n\in\Z} I_{n,0} = \text{supp}(u_0) \quad\text{and}\quad I_\pm = \bigcup_{n\in\Z} I_{n,\pm} = \text{supp}(u_\pm) .\]
%\end{definition}

\section{The Boussinesq model on the necklace graph}

\label{sec3}

%So far, very few evolutionary problems defined by nonlinear PDEs posed on quantum graphs have been considered.
%Examples are the NLS equation, the Burgers equation, or the cubic Klein-Gordon equation, cf. \cite{GPS16,GSU20}. In particular, almost no results on the qualitative behavior of solutions to the initial value problem can be found in the literature.
%To this rather short list of examples we would like to add another example in this paper.
As already said we are interested in the qualitative behavior of the solutions of  the Boussinesq equation 
\begin{equation}\label{Bouss}
\partial_t^2 u = \partial_x^2 u + \partial_t^2 \partial_x^2 u + \partial_x^2 (u^2),
\end{equation}
posed  on the infinite periodic necklace  graph $\Gamma$ 
sketched in Figure \ref{fig2}. 
This equation is a phenomenological model of the so-called water wave problem, cf. \cite{Cr85}.

First of all, the Boussinesq equation on the infinite periodic necklace graph has to be defined.
\begin{definition}
The Boussinesq equation posed  on the infinite periodic necklace  graph $\Gamma$  is 
the vector-valued PDE 
\begin{equation}\label{eq. Bouss1}
\partial_t^2 U = \partial_x^2 U + \partial_t^2 \partial_x^2 U + \partial_x^2 (U^2),\quad t\in\R, \; x\in\R\setminus (\pi\Z),
\end{equation}
for  $U=(u_0,u_+,u_-)$ 
with the Kirchhoff boundary conditions \eqref{eq. KH1}--\eqref{eq. KH2} at the vertices  $x\in\pi\Z$ and where the quadratic nonlinearity is defined component-wise, i.e., $U^2=(u_0^2,u_+^2,u_-^2)$.
\end{definition}

For  proving a local existence and uniqueness result we need a few preparations.
\begin{definition}
Our basic space $\mathcal{L}^2$ is defined by
\[ \mathcal{L}^2 = \{ U = (u_0,u_+,u_-)\in (L^2(\R))^3 :\; \text{supp}(u_{n,j})=I_{n,j},\;n\in\Z,\; j\in\{0,+,-\} \} \]
equipped with the scalar product
\[ \langle U , V \rangle_{\mathcal{L}^2} = \sum_{n\in\Z} (( u_{n,0} , v_{n,0} )_{L^2(I_{n,0})} + ( u_{n,+} , v_{n,+} )_{L^2(I_{n,+})} + ( u_{n,-} , v_{n,-} )_{L^2(I_{n,-})}). \]
We further define the space $\mathcal{H}^m$ with $m\in\N$ by
\[ \mathcal{H}^m := \{ U \in \mathcal{L}^2 :\; u_{n,j} \in H^m(I_{n,j}),\; n\in\Z,\; j\in\{0,+,-\},\; \eqref{eq. KH1},\eqref{eq. KH2} \;\text{satisfied} \}, \]
equipped with the norm
\[ \| U\|_{\mathcal{H}^m} = \left( \sum_{n\in\Z} (\| u_{n,0} \|_{H^m(I_{n,0})}^2 + \| u_{n,+} \|_{H^m(I_{n,+})}^2 + \| u_{n,-} \|_{H^m(I_{n,-})}^2 )\right)^\frac{1}{2}. \]
For $ \mathcal{H}^1 $ we only assume the condition \eqref{eq. KH1}.
Moreover, we introduce the norm 
$$
\| U\|_{L^{\infty}} =  \sum_{n\in\Z} (\| u_{n,0} \|_{L^{\infty}} + \| u_{n,+} \|_{L^{\infty}} + \| u_{n,-} \|_{L^{\infty}}).
$$
\end{definition}
We recall from \cite[Section 3]{GPS16}:
\begin{lemma}
The domain of definition of the operator $\mathcal{A}^2=-\partial_x^2 $ in $\mathcal{L}^2$ is given by $\mathcal{H}^2$.
\end{lemma}
%and: %We recall \cite[Lemma 3.2]{GPS16}:
\begin{lemma}\label{Lem. -dx^2}
The operator $\mathcal{A}^2 : \mathcal{H}^2\rightarrow \mathcal{L}^2$ is self-adjoint and positive semidefinite with respect to the $\mathcal{L}^2$-scalar product.
There exists a self-adjoint and positive semidefinite root $ \mathcal{A} $.
\end{lemma}
Direct consequences of these facts are 
\begin{lemma}\label{Lem. -dx^21}
The operator $(I+\mathcal{A}^2 )^{-1}: \mathcal{L}^2\rightarrow \mathcal{H}^2$ is bounded.
\end{lemma}
%Therefore, and since $ (I+\mathcal{A}^2 )^{-1} $ and $ \mathcal{A} $ commutate, we have 
\begin{lemma}\label{Lem. -dx^22}
The operator $\mathcal{B}^2 = (I+\mathcal{A}^2 )^{-1} \mathcal{A}^2 : \mathcal{H}^2\rightarrow \mathcal{H}^2$ is bounded.
\end{lemma}
Hence, we can rewrite \eqref{eq. Bouss1} as
\begin{equation}\label{eq. Bouss2}
\partial_t^2 U = -\mathcal{B}^2 U - \mathcal{B}^2  (U^2).
\end{equation}
In order to write this equation as first order system we use
\begin{corollary}\label{Lem. L^2}
The operator $\mathcal{B}^2 $ is self-adjoint and positive semidefinite with respect to the $\mathcal{L}^2$-scalar product.
There exists a unique self-adjoint and positive semidefinite root $\mathcal{B}$ of $\mathcal{B}^2$.  The operator $ \mathcal{B} $ is bounded in $ \mathcal{H}^2 $.
\end{corollary}
\begin{proof}
By Lemma \ref{Lem. -dx^2}, $I+\mathcal{A}^2 $ is self-adjoint and positive semidefinite, 
too, and so is its inverse $(I+\mathcal{A}^2 )^{-1}$. This is also true for  $\mathcal{B}^2 = (I+\mathcal{A}^2)^{-1}\mathcal{A}^2$, because of $ (I+\mathcal{A}^2)^{-1}\mathcal{A}^2 = \mathcal{A}^2(I+\mathcal{A}^2)^{-1} $.
By construction $ \mathcal{B} $ is bounded in $ \mathcal{L}^2 $.
Since obviously 
$ 1 + \mathcal{A}^2 $ and $ \mathcal{B} $ commute,
we have 
$$ 
\| \mathcal{B} u \|_{\mathcal{H}^2} \leq C \|(1 + \mathcal{A}^2) \mathcal{B} u \|_{\mathcal{L}^2} \leq C \| \mathcal{B}  (1 + \mathcal{A}^2) u \|_{\mathcal{L}^2}
\leq C \| (1 + \mathcal{A}^2) u \|_{\mathcal{L}^2}
\leq  C \|  u \|_{\mathcal{H}^2}
$$
and so the  boundedness of $ \mathcal{B} $  in $ \mathcal{H}^2 $.

%\[ (I+\mathcal{A}^2)^{-1}\mathcal{A}^2= (\mathcal{A}^{-2}(I+\mathcal{A}^2))^{-1} = (\mathcal{A}^{-2}+I)^{-1} = ((I+\mathcal{A}^2)\mathcal{A}^{-2})^{-1} = \mathcal{A}^2(I+\mathcal{A}^2)^{-1} \]
%commutate. 
\end{proof}

The validity of the KdV approximation  will be proved in the space $ \mathcal{H}^1 $.
Obviously the space $ \mathcal{H}^1 $ is closed under pointwise multiplication, i.e.,
there is a $ C > 0 $ such that for all $ u,v \in  \mathcal{H}^1 $ we 
have 
$$ 
\| u v \|_{\mathcal{H}^1} \leq C \| u  \|_{\mathcal{H}^1} \|  v \|_{\mathcal{H}^1} .
$$ 
Moreover, we have 
$$ 
\| \partial_x u \|^2_{\mathcal{L}^2} =  \langle \partial_x u, \partial_x u \rangle_{\mathcal{L}^2} = -  \langle  u, \partial_x^2 u \rangle_{\mathcal{L}^2} 
=   \langle  u, \mathcal{A}^2 u \rangle_{\mathcal{L}^2} 
=   \langle   \mathcal{A} u, \mathcal{A} u \rangle_{\mathcal{L}^2} 
= \|  \mathcal{A} u \|^2_{\mathcal{L}^2},
$$ 
where we used integration by parts at the second step. 
No boundary terms occur at the vertex points due to the boundary conditions 
\eqref{eq. KH1} and \eqref{eq. KH2}. 
In the third step we used the definition 
of the operator $ \mathcal{A}^2 $. 
In the fourth step we used the self-adjointness 
of the operator $ \mathcal{A} $. 
Therefore, the $ \mathcal{H}^1 $-norm of a function $ u $ can be estimated by 
$ \|u\|_{\mathcal{L}^2} + \|\mathcal{A}  u \|_{\mathcal{L}^2} $ and vice versa.  This implies 
\begin{eqnarray*}
\| \mathcal{B} u \|_{\mathcal{H}^1} & \leq & C 
(\|  \mathcal{B} u \|_{\mathcal{L}^2} + 
\| \mathcal{A} \mathcal{B} u \|_{\mathcal{L}^2}) \\& \leq &
C (\|  \mathcal{B} u \|_{\mathcal{L}^2} + 
\| \mathcal{B} \mathcal{A}  u \|_{\mathcal{L}^2})   \leq 
C (\|  u \|_{\mathcal{L}^2} + 
\| \mathcal{A}  u \|_{\mathcal{L}^2}) \leq 
C \|  u \|_{\mathcal{H}^1}
\end{eqnarray*}
and so the  boundedness of $ \mathcal{B} $  in $ \mathcal{H}^1 $.

Hence, with Corollary \ref{Lem. L^2}
 we transform \eqref{eq. Bouss2} into a first order system
\begin{equation}\label{eq. Bouss3}
\partial_t W = \Lambda W + F(W),
\end{equation}
where $W=(U,V)$,
\begin{align*}
\Lambda = \left(\begin{matrix} 0 & i \mathcal{B} \\ i \mathcal{B} & 0 \end{matrix}\right) \quad \text{and} \quad F(W) = \left(\begin{array}{c} 0 \\ i\mathcal{B}(U^2) \end{array}\right).
\end{align*}
For proving the local existence and uniqueness of solutions we use   the Picard-Lindel\"of theorem. 
Since the spaces $\mathcal{H}^s$ for $ s =1,2 $ is closed under pointwise multiplication, cf.
\cite[Lemma 3.1]{GPS16} for $\mathcal{H}^2$ and  for $\mathcal{H}^1$ see above,
 the right hand side of \eqref{eq. Bouss3} is locally 
Lipschitz-continuous from $ \mathcal{H}^s\rightarrow \mathcal{H}^s$ 
for $ s =1,2 $
and so we have
\begin{theorem}
Let $ s =1,2 $.
For each initial condition $ W_0 \in (\mathcal{H}^s)^2 $
there exists a time $T_0 = T_0(\|W_0\|_{(\mathcal{H}^s)^2}) > 0$ and a unique solution $W \in  C([-T_0,T_0],(\mathcal{H}^s)^2)$ of \eqref{eq. Bouss3} with $ W|_{t=0} = W_0 $.
\end{theorem}

\begin{remark}{\rm
The local  existence and uniqueness of solutions combined with the 
subsequent error estimates 
yields the existence and uniqueness of solutions for all $ t \in [0,T_0/\varepsilon^3] $. 
}
\end{remark}

\section{Floquet-Bloch spectrum}\label{subsecspec}

\label{sec4}

The KdV equation describes non-oscillatory Fourier or Bloch modes 
at the Fourier or Bloch wave number $ k = 0 $ or $ l = 0 $.
Therefore, we analyze the linearized system first. 
In the following we restrict ourselves to the invariant subspace of symmetric 
functions, i.e., to functions satisfying 
\begin{equation} \label{usym}
u_+ = u_-.
\end{equation}
For such functions the necklace graph can be identified with the real line
with boundary conditions at $ \{x = n \pi : n \in \Z \} $ coming from 
\eqref{eq. KH1} and \eqref{eq. KH2}.
We consider  the linear part 
\begin{equation}\label{eq. Bouss1lin}
\partial_t^2 U = \partial_x^2 U + \partial_t^2 \partial_x^2 U ,\qquad t\in\R, \; x\in\R\setminus (\pi\Z),
\end{equation}
of \eqref{eq. Bouss1},
respectively  the linearization of \eqref{eq. Bouss2}, namely
\begin{equation}\label{eq. Bouss2lin}
\partial_t^2 U = -\mathcal{B}^2 U =  -(I+\mathcal{A}^2 )^{-1} \mathcal{A}^2 U.
\end{equation}
Inserting the ansatz
\[ U(x,t)=W_B(x) e^{i\omega t} \]
into \eqref{eq. Bouss2lin} yields the spectral problem
\[ \mathcal{B}^2 W_B = (I+\mathcal{A}^2 )^{-1} \mathcal{A}^2  W_B = \omega^2 W_B = \mu W_B. \]
The components of $W_B=(w_0,w_+,w_-)$ satisfy the Kirchhoff boundary conditions \eqref{eq. KH1}--\eqref{eq. KH2} and have their supports in $(I_0,I_+,I_-)$. 
The eigenfunctions $W_B$ of $\mathcal{B}^2$ can be written in the form of Bloch waves 
\[ W_B(x)=e^{ilx}f(l,x),\quad l,x\in\R. \]
The functions $f=(f_0,f_+,f_-)$ satisfy 
\[ f(l,x)=f(l,x+2\pi),\quad f(l,x)=f(l+1,x)e^{ix},\quad l,x\in\R. \]
Therefore, we can restrict the definition of $f(l,x)$ to $x\in\T_{2\pi}=\R /(2\pi\Z)$ and $l\in\T_{1}=\R / \Z $.
We now consider the eigenvalue problem
\begin{equation}\label{eq. EWP*}
\widetilde{\mathcal{B}}_l^2 f = \mu(l) f,\quad x\in\T_{2\pi},
\end{equation}
for fixed $l\in\T_{1}$ with
\begin{equation}\label{eq. B tilde}
\widetilde{\mathcal{B}}^2_l=-(1-(\partial_x + il)^2)^{-1}(\partial_x + il)^2
\end{equation}
subject to the boundary conditions
\begin{equation}\label{eq. KH1*}
\begin{cases}
f_0(l,\pi)=f_+(l,\pi)=f_-(l,\pi), \\
f_0(l,2\pi)=f_+(l,2\pi)=f_-(l,2\pi)
\end{cases}
\end{equation}
and
\begin{equation}\label{eq. KH2*}
\begin{cases}
(\partial_x + il)f_0(l,\pi)=(\partial_x + il)f_+(l,\pi)+(\partial_x + il)f_-(l,\pi), \\
(\partial_x + il)f_0(l,2\pi)=(\partial_x + il)f_+(l,2\pi)+(\partial_x + il)f_-(l,2\pi),
\end{cases}
\end{equation}
which can be derived from \eqref{eq. KH1}--\eqref{eq. KH2} using the $2\pi$-periodicity of $f(l,\cdot)$. The functions $f_0(l,\cdot)$ and $f_\pm (l,\cdot)$ have their supports in $I_{0,0}\subset\T_{2\pi}$ and $I_{0,\pm}\subset\T_{2\pi}$, respectively.

For solving the eigenvalue problem \eqref{eq. Bouss2lin} we rewrite it as 
$$ 
 \mathcal{A}^2  W_A = \omega^2 (I+\mathcal{A}^2 )  W_A
$$ 
and use that the solution of the eigenvalue problem  $  \mathcal{A}^2  W_A = \lambda W_A $ has already been solved, 
cf. \cite{GPS16}.
The eigenfunctions $W_A$ of $\mathcal{A}^2$ can be written in the form of Bloch waves 
\[ W_A(x)=e^{ilx}f(l,x),\quad l,x\in\R, \]
with $ f $ having the  properties as above.
For fixed $ l \in \T_{1} $ we find 
\begin{equation}\label{eq. EVP A}
\widetilde{\mathcal{A}}_l^2 f(l,\cdot) = \lambda(l) f(l,\cdot) = \mu(l) (I+\widetilde{\mathcal{A}}_l^2 )   f(l,\cdot) = \mu(l) (1 +  \lambda(l)) f(l,\cdot),
\end{equation}
where 
$ \widetilde{\mathcal{A}}_l^2= -(\partial_x + il)^2 $.
Therefore,
\begin{equation}\label{eq. relation1}
\mu(l) = \frac{\lambda(l)}{1+\lambda(l)}.
\end{equation}
The curves $ l \to \lambda(l) $  have been computed for instance in \cite[Section 2]{Pelinovsky*}.
In the symmetric case \eqref{usym},  they correspond to the real roots $ \rho_{1,2} $ of  $ \rho^2 - {\rm tr}(M)(\lambda) \rho +1 = 0  $ with $\rho=e^{2 \pi il}$. 
 Here 
\begin{equation}
\label{transcendental}
{\rm tr}(M)(\lambda) := 
2 \cos(\pi \sqrt{\lambda})^2 - (5/2) \sin(\pi \sqrt{\lambda})^2
%\frac{1}{4} \left[ 9 \cos(2 \pi \sqrt{\lambda}) - 1\right].
\end{equation}
 is the trace of the monodromy matrix $M$ associated with the eigenvalue problem $  \mathcal{A}^2  W_A = \lambda W_A $.  
Real roots are obtained if $ {\rm tr}(M)(\lambda) \in [-2,2] $.  

\begin{figure}[ht]
\centering
\includegraphics[scale=0.42]{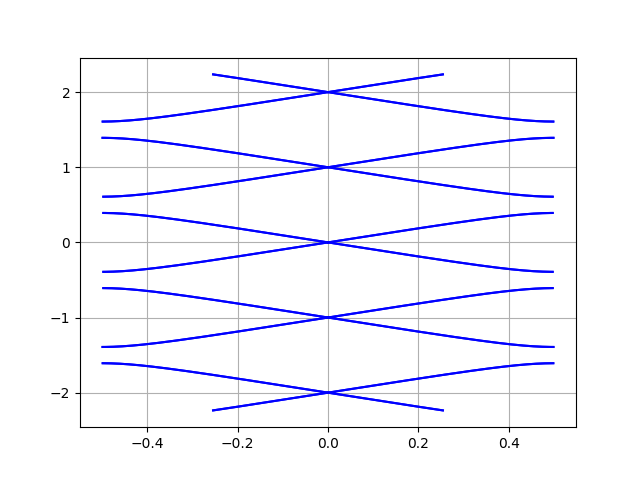} \qquad
\includegraphics[scale=0.42]{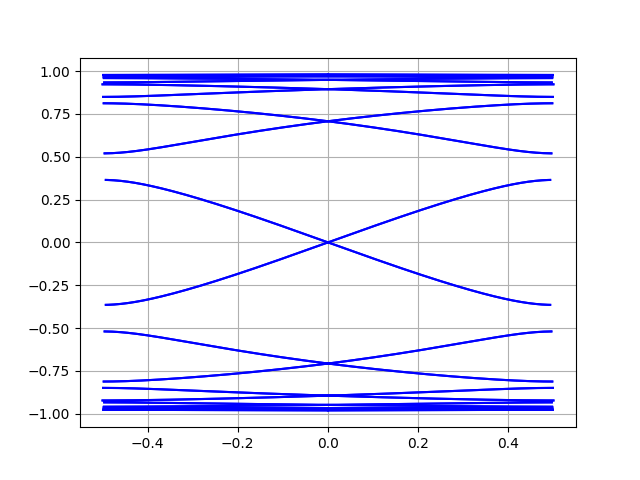}
\caption{The spectral bands 
$ l \mapsto \pm\sqrt{\lambda(l)} $ to the spectral problem \eqref{eq. EVP A} and the spectral bands 
$l \mapsto  \pm\omega = \pm\sqrt{\frac{\lambda(l)}{1+\lambda(l)}} $ to the spectral problem \eqref{eq. EWP*} in the interval $l\in [-1/2,1/2]$ in the symmetric situation.}
\label{Fig. SpektrumPlot}
\end{figure}

Hence, due to \eqref{eq. relation1}, we can also compute the eigenvalues $\mu(l)$ of the operator $\widetilde{\mathcal{B}}_l^2$. The curves $ l \mapsto \omega(l) = \pm\sqrt{\frac{\lambda(l)}{1+\lambda(l)}} $ 
are plotted in Figure \ref{Fig. SpektrumPlot}. We see that two curves 
go through the origin. The associated modes  can be described by KdV equations, as
we will see in the following.

In order to properly define the eigenfunctions to the eigenvalue problem \eqref{eq. EWP*} we define 

\begin{definition}
For fixed $l\in\T_{1}$ let 
\[ L_\Gamma^2 := \{ \widetilde{U} = (\widetilde{u}_0,\widetilde{u}_+,\widetilde{u}_-)\in (L^2(\T_{2\pi}))^3 : \; \text{supp}(\widetilde{u}_j)=I_{0,j},\; j\in \{ 0,+,- \} \} \]
with the scalar product
\[ \langle \widetilde{U},\widetilde{V} \rangle_{L_\Gamma^2} = ( \widetilde{u}_{0} , \widetilde{v}_{0} )_{L^2(I_{0,0})} + ( \widetilde{u}_{+} , \widetilde{v}_+ )_{L^2(I_{0,+})} + ( \widetilde{u}_{-} , \widetilde{v}_{-} )_{L^2(I_{0,-})} .\]
\end{definition}

%and
%\[ H_\Gamma^2 (l) := \{ \widetilde{U}\in L_\Gamma^2 : \; \widetilde{u}_j \in H^2(I_{0,j}),\; j\in \{ 0,+,- \},\; \eqref{eq. KH1*}-\eqref{eq. KH2*} \; \text{are satisfied} \} \]
%equipped with the norm
%\[ \| \widetilde{U} \|_{H_\Gamma^2 (l)} = \left( \| \widetilde{u}_0 \|_{H^2 (I_{0,0})}^2 + \| \widetilde{u}_+ \|_{H^2 (I_{0,+})}^2 + \| \widetilde{u}_- \|_{H^2 (I_{0,-})}^2 \right)^\frac{1}{2}. \]
%\end{definition}

From \cite[Lemma 2.2]{GPS16} we know that the operator $\widetilde{\mathcal{A}}_l^2$ is self-adjoint and positive semidefinite in $L_\Gamma^2$ for fixed $l\in\T_1$. By the same reasoning as in the proof of Corollary \ref{Lem. L^2}, this is also true for $\widetilde{\mathcal{B}}^2_l$. Because of the spectral theorem for self-adjoint operators, cf. \cite[Chapter 7]{Reed}, for fixed $l\in\T_1$ there exists a countable orthonormal basis of $L_\Gamma^2$ of eigenfunctions
\begin{align*}
f^{(m)}(l,x)= \left( f^{(m)}_0(l,x), f^{(m)}_+(l,x), f^{(m)}_-(l,x) \right)
\end{align*}
of $\widetilde{\mathcal{B}}^2_l$ with associated nonnegative eigenvalues $\mu^{(m)}(l)$ satisfying $\mu^{(m)}(l)\leq \mu^{(m+1)}(l)$ for all $m \in \N$.

For the derivation of the KdV equation we need some additional properties which 
due to \eqref{eq. EVP A}--\eqref{eq. relation1}  are direct corollaries of \cite[Lemma 3.1, Lemma 3.2]{BDS19}.
\begin{lemma}\label{Lem. Eig. L}
For $l=0$, the operator $\widetilde{\mathcal{B}}_0^2$ has the simple eigenvalue $\mu^{(1)}(0)=0$ to the eigenfunction $f^{(1)}(0,x)=(3 \pi)^{-1/2}\, 1_3$, where $1_3=(1,1,1)^T$.
\end{lemma}

 It is well known  that the curves $l\mapsto \mu^{(m)}(l)$ and $l\mapsto f^{(m)}(l,\cdot)$ are smooth with respect to $l$ for simple eigenvalues. Thus there exists a $\delta_0 > 0$ such that the eigenvalue $\mu^{(1)}(l)$ is separated from the rest of the spectrum for all $l\in [ -\delta_0,\delta_0 ]$.

\begin{lemma}\label{Lem. Eigfkt}
The curve $l\mapsto \mu^{(1)} (l),\; l\in [-\delta_0,\delta_0]$ is an even real-valued function. The associated curve of eigenfunctions satisfies $f^{(1)}(l,x) = \overline{f^{(1)}(-l,x)}$ and possesses an expansion
\[ f^{(1)}(l,x) = \sum_{j=0}^\infty (il)^j g_j(x) \]
with $g_0(x)=(3 \pi)^{-1/2}\,1_3$, $\langle g_j(\cdot),1_{3} \rangle_{L_\Gamma^2} = 0$ for $j\geq 1$ and
\[ g_{2j}(x)=g_{2j}(-x)\in\R^3,\quad g_{2j+1}(x)=-g_{2j+1}(-x)\in\R^3. \]
%With the evolution of $f^{(1)}$ obviously holds
%\begin{equation}\label{eq. Formel f(1)}
%f^{(1)}(l-m)f^{(1)}(m) = f^{(1)}(l) + \mathcal{O}(l^2+(l-m)^2+m^2).
%\end{equation}

\end{lemma}

\section{The Boussinesq equation in Bloch space}

\label{sec5}

The derivation of amplitude, envelope, or modulation equations 
for spatially periodic systems 
heavily relies
on an expansion of the original system in  Bloch space. 
We start this section by recalling some basic properties of  
Bloch transform, cf. for instance  \cite[Chapter 3.2]{Gallay} or \cite[Chapter 2.1.2]{Pelinovsky}.

\begin{definition}
Bloch transform $\mathcal{T}$ is defined by
\[ \widetilde{u}(l,x) = (\mathcal{T}u)(l,x) = \sum_{n\in\Z} u(x+2\pi n)e^{-ilx-2\pi inl} \]
and its inverse by
\[ u(x) = (\mathcal{T}^{-1}\widetilde{u})(x) = \int_{-\frac{1}{2}}^\frac{1}{2} e^{ilx}  \widetilde{u}(l,x)\; dl. \]
with  $ \widetilde{u} $ satisfying the continuity conditions
\[ \widetilde{u}(l,x) = \widetilde{u}(l,x+2\pi),\quad \widetilde{u}(l,x) = \widetilde{u}(l+1,x)e^{ix}. \]
\end{definition}
Thus  $\mathcal{T}$ can be extended from $(l,x)\in\T_1\times\T_{2\pi}$ to $(l,x)\in\R\times\R$.
Its properties are summarized in the following lemma.
\begin{lemma}
a)  For real-valued $u$ it additionally holds
\[ \widetilde{u}(l,x) = \bar{\widetilde{u}}(-l,x). \]

b) Bloch transform $\mathcal{T}$ is an isomorphism between $H^s(\R)$ and $L^2(\T_1,H^s(\T_{2\pi}))$, i.e., we have
\[ \| u \|_{H^s(\R)} = \left( \int_{-\frac{1}{2}}^\frac{1}{2} \| \widetilde{u}(l,\cdot) \|_{H^s(\T_{2\pi})}^2 \; dl\right)^\frac{1}{2}. \]

c) Furthermore,
\[ \widetilde{\partial_x u} (l,x) = (\partial_x + il)\widetilde{u}(l,x). \]

d) Bloch transform and  Fourier transform are related through
\[ \widetilde{u}(l,x) = \sum_{j\in\Z} e^{ijx} \hat{u}(l+j). \]

e) The Bloch transform of a product of two functions $u$ and $v$ is the convolution of the respective Bloch transformations $\widetilde{u}$ and $\widetilde{v}$ with respect to 
the Bloch wave number $l$, i.e.
\[ (\widetilde{uv})(l,x) = (\widetilde{u}\ast\widetilde{v})(l,x) = \int_{-1/2}^{1/2} \widetilde{u}(l-l',x)\widetilde{v}(l',x) dl'. \]

f) Finally, for a $2\pi$-periodic function $\chi:\R\rightarrow\R$ we have
\[ \mathcal{T}(\chi u)(l,x) = \chi(x)(\mathcal{T}u)(l,x). \]
\end{lemma}
With the help of these properties, we can transfer  equation \eqref{eq. Bouss1}, resp. \eqref{eq. Bouss2}, into Bloch space by applying the Bloch transform $\mathcal{T}$.
We obtain
\begin{equation}\label{eq. Bouss1*}
\partial_t^2 \widetilde{U}(t,l,x) = -\widetilde{\mathcal{B}}_l^2\widetilde{U}(t,l,x)-\widetilde{\mathcal{B}}_l^2N(\widetilde{U})(t,l,x),
\end{equation}
where $\widetilde{\mathcal{B}}_l^2$ is given by \eqref{eq. B tilde} and the convolution in the nonlinearity  is applied component-wise, i.e. 
\[ N(\widetilde{U})(t,l,x) = (\widetilde{U}\ast\widetilde{U})(t,l,x)=(\widetilde{u}_0\ast\widetilde{u}_0,\widetilde{u}_+\ast\widetilde{u}_+,\widetilde{u}_-\ast\widetilde{u}_-)(t,l,x). \]
The functions $\widetilde{U}(t,l,x)=(\widetilde{u}_0,\widetilde{u}_+,\widetilde{u}_-)$ satisfy
\begin{equation} 
\label{period1}
\widetilde{U}(t,l,x)= \widetilde{U}(t,l,x+2\pi),\quad \widetilde{U}(t,l,x) = \widetilde{U}(t,l+1,x)e^{ix}. \end{equation}
%\begin{remark}{\rm 
As in \cite[Subsection 4.2]{GPS16} we define periodic truncation functions
\begin{align*}
\chi_j (x) = 
\begin{cases}
1, & x\in I_j \\
0, & \text{else}
\end{cases}, \quad j\in \{ 0,+,- \}.
\end{align*}
We  obviously have
\[ \mathcal{T}(u_j)(l,x) = \mathcal{T}(\chi_j u_j)(l,x) = \chi_j (x) \mathcal{T}(u_j)(l,x) ,\]
and thus the support of $\mathcal{T}(u_j)(l,\cdot)$ is contained in $I_j$  for each $j\in \{ 0,+,- \}$.
%}\end{remark}

Now that we are no longer in physical space but in Bloch space we need a replacement 
for the space $\mathcal{H}^2$.
\begin{definition} 
%For fixed $l\in\T_{1}$ we define 
%\[ L_\Gamma^2 := \{ \widetilde{U} = (\widetilde{u}_0,\widetilde{u}_+,\widetilde{u}_-)\in (L^2(\T_{2\pi}))^3 : \; \text{supp}(\widetilde{u}_j)=I_{0,j},\; j\in \{ 0,+,- \} \} \]
%with the scalar product
%\[ \langle \widetilde{U},\widetilde{V} \rangle_{L_\Gamma^2} = ( \widetilde{u}_{0} , \widetilde{v}_{0} )_{L^2(I_{0,0})} + ( \widetilde{u}_{+} , \widetilde{v}_+ )_{L^2(I_{0,+})} + ( \widetilde{u}_{-} , \widetilde{v}_{-} )_{L^2(I_{0,-})} .\]
%Then 
For $ s = 1,2 $ we define the space
\[ \widetilde{\mathcal{H}}^s := \{ \widetilde{U}\in L^2(\T_1,L_\Gamma^2): \; \widetilde{u}_j \in L^2(\T_1,H^2(I_{0,j})),\; j\in\{ 0,+,- \},\; \eqref{eq. KH1*}-\eqref{eq. KH2*}\;\text{satisfied} \} \]
equipped with the norm
\[ \| \widetilde{U} \|_{\widetilde{\mathcal{H}}^s} = \left(\int_{-\frac{1}{2}}^\frac{1}{2} \left( \| \widetilde{u}_0(l,\cdot) \|_{H^s(I_{0,0})}^2 + \| \widetilde{u}_+(l,\cdot) \|_{H^s(I_{0,+})}^2 + \| \widetilde{u}_-(l,\cdot) \|_{H^s(I_{0,-})}^2 \right)\; dl \right)^\frac{1}{2}, \]
where for $ \widetilde{\mathcal{H}}^1 $ only \eqref{eq. KH1*} should hold.
\end{definition}
%Obviously $\widetilde{\mathcal{B}}_l^2$ maps the space $\widetilde{\mathcal{H}}^2$ into itself.

With \cite[Lemma 4.2]{GPS16} we recall another property of Bloch transform.
\begin{lemma}\label{Lem. Isom}
For $ s = 1,2 $ the Bloch transform $\mathcal{T}$ is an isomorphism between $\mathcal{H}^s$ and $\widetilde{\mathcal{H}}^s$.
\end{lemma}

\section{Derivation of the KdV equation}

\label{sec6}

For the  derivation of  the KdV equation we divide the Boussinesq equation \eqref{eq. Bouss1*} in  Bloch space in two parts, cf. \cite[Chapter 3]{BDS19}, namely 
\begin{equation}\label{eq. Ansatz2}
\widetilde{U}(t,l,x) = \chi_{[-\frac{1}{2}\delta_0,\frac{1}{2}\delta_0]}(l) \widetilde{V}(t,l)f^{(1)}(l,x)+ \widetilde{U}^\perp (t,l,x)
\end{equation}
with 
the orthogonality condition
\[ \langle f^{(1)}(l,\cdot),\widetilde{U}^\perp (t,l,\cdot) \rangle_{L_\Gamma^2} = 0\]
for 
$  l\in \left[-\frac{1}{2}\delta_0,\frac{1}{2}\delta_0 \right] $ and $ \delta_0 > 0 $ from Lemma \ref{Lem. Eigfkt}, leading 
to the uniqueness of the decomposition. We find
\begin{equation}\label{eq. Bouss2*}
\begin{array}{rcl}
\partial_t^2 \widetilde{V}(t,l) &=& -\mu^{(1)}(l)\widetilde{V}(t,l)- P^c_l\widetilde{\mathcal{B}}_l^2 N(\widetilde{V},\widetilde{U}^\perp)(t,l), \\[1mm]
\partial_t^2 \widetilde{U}^\perp(t,l,x) &=& -\widetilde{\mathcal{B}}_l^2  \widetilde{U}^\perp(t,l,x) - P^s_l\widetilde{\mathcal{B}}_l^2 N(\widetilde{V},\widetilde{U}^\perp)(t,l,x),
\end{array}
\end{equation}
with $N(\widetilde{V},\widetilde{U}^\perp)=N(\widetilde{U})$, where
\[ P^c_l \widetilde{U}(t,l) = % \frac{1}{3\pi}
\chi_{[-\frac{1}{2}\delta_0,\frac{1}{2}\delta_0]}(l) \langle f^{(1)}(l,\cdot),\widetilde{U}(t,l,\cdot)\rangle_{L_\Gamma^2} \]
and
\[ P^s_l \widetilde{U}(t,l,x) =  \widetilde{U}(t,l,x) - P^c_{l}  \widetilde{U}(t,l) f^{(1)}(l,x). \]
The nonlinearity $P^c_{l}  \widetilde{\mathcal{B}}_l^2 N(\widetilde{V},\widetilde{U}^\perp)$ can be written as
\begin{align*}
& P^c_l \widetilde{\mathcal{B}}_l^2 N(\widetilde{V},\widetilde{U}^\perp)(t,l)  \\
&\qquad\qquad=  \int_{\T_1} \beta(l,l-m,m) \widetilde{V}(t,l-m)\widetilde{V}(t,m) \; dm + \widetilde{\mathcal{B}}_l^2 N_\text{rest}(\widetilde{V},\widetilde{U}^\perp)(t,l),
\end{align*}
where
\begin{equation}\label{eq. beta}
\beta(l,l-l',l') := %\frac{1}{3\pi} w^{(1)}(l)
 \langle f^{(1)}(l,\cdot),\widetilde{\mathcal{B}}_l^2 (f^{(1)}(l-l',\cdot)f^{(1)}(l',\cdot))\rangle_{L_\Gamma^2}.
\end{equation} 
The term $\widetilde{\mathcal{B}}_l^2 N_\text{rest}(\widetilde{V},\widetilde{U}^\perp)$ contains all remaining convolutions, each of which has at least one $\widetilde{U}^\perp$ as a factor, so obviously 
$$\widetilde{\mathcal{B}}_l^2 N_\text{rest}(\widetilde{V},0)(t,l)=0.$$ 
Furthermore, we have the following lemma. 
\begin{lemma}\label{Lem. beta}
We have
\[ \left| \beta (l,l-l',l') -   \frac{1}{2}\partial_l^2 \mu^{(1)}(0)\, l^2 \right| = \mathcal{O}(l^2(l^2+(l-l')^2+l'^2)). \]
%with $ \nu_2 $ a fixed constant.
\end{lemma}
\begin{proof}
The result follows after a straightforward calculation from Lemma \ref{Lem. Eigfkt}, analogously  as in the proof of  \cite[Lemma 3.3]{BDS19}.
%
%
%
%We have
%\[ f^{(1)}(l,x) = 1_3+ilg_1(x) + \mathcal{O}(l^2) \]
%with $g_1(x)\in\R^3$. That delivers
%\begin{align*}
%\beta(l,l-m,m) &= \frac{1}{3\pi} w^{(1)}(l) \left\langle f^{(1)}(l,\cdot),f^{(1)}(l-m,\cdot)f^{(1)}(m,\cdot)\right\rangle_{L_\Gamma^2} \\
%&= \frac{1}{3\pi} w^{(1)}(l)\sum_{j\in\{ 0,+,- \}} \int _{I_{0,j}} \overline{f^{(1)}_j(l,x)} f^{(1)}_j(l-m,x)f^{(1)}_j(m,x) \; dx \\
%&= \frac{1}{3\pi}w^{(1)}(l) \sum_{j\in\{ 0,+,- \}} \int _{I_{0,j}} (1-ilg_{1,j}(x) + \mathcal{O}(l^2)) \\
%&\qquad\times (1+i(l-m)g_{1,j}(x) + \mathcal{O}((l-m)^2)) (1+img_{1,j}(x) + \mathcal{O}(m^2)) \; dx \\
%&\qquad+ \mathcal{O}(l^2+(l-m)^2+m^2) \; dx \\
%&= \frac{1}{3\pi} w^{(1)}(l) \left(\sum_{j\in\{ 0,+,- \}} \int _{I_{0,j}} 1 \; dx + \mathcal{O}(l^2+(l-m)^2+m^2)\right) \\
%&= w^{(1)}(l) \left(1 + \mathcal{O}(l^2+(l-m)^2+m^2)\right) \\
%&= \frac{1}{2}\partial_l^2 w^{(1)}(0)l^2 + \mathcal{O}(l^2(l^2+(l-m)^2+m^2))
%\end{align*}
\end{proof}
For the derivation of the KdV equation  we make the ansatz
\begin{equation}\label{eq. Ansatz3  V}
\widetilde{V}_\text{KdV}(t,l) = \varepsilon \widetilde{A} (T,K)\,{\bf E}(t,l),
\end{equation}
where $T=\varepsilon^3 t$, $l=\varepsilon K$, and ${\bf E}(t,l)=e^{-ilct}$. Substituting the ansatz into the first equation of \eqref{eq. Bouss2*} with $\widetilde{U}^\perp_\text{KdV}(t,l,x) = 0$, results in
\begin{align*}
&\varepsilon^7 \partial_T^2 \widetilde{A} \,{\bf E} - 2i\varepsilon^5 cK \partial_T \widetilde{A}\,{\bf E} - \varepsilon^3 c^2 K^2 \widetilde{A}\,{\bf E} \\
&\qquad = -\varepsilon \widetilde{A}\,{\bf E} \left(\frac{1}{2} \partial_l^2 \mu^{(1)}(0) \varepsilon^2 K^2 + \frac{1}{24} \partial_l^4 \mu^{(1)}(0)\varepsilon^4 K^4 + \mathcal{O}(\varepsilon^6)\right) \\
&\qquad\qquad - \varepsilon^5 \frac{1}{2} \partial_l^2 \mu^{(1)}(0) \int_{-\frac{1}{2\varepsilon}}^\frac{1}{2\varepsilon} K^2 \widetilde{A}(T,K-M)\widetilde{A}(T,M)\; dM\, {\bf E},
\end{align*}
where we expanded $\mu^{(1)}$ around $l=\varepsilon K =0$ and made the substitution $m =\varepsilon M $ in the integral above \eqref{eq. beta}. Equating the coefficient in front  of $\varepsilon^3 {\bf E}$ yields
\begin{equation}\label{eq. c}
c^2 = \frac{1}{2}\partial_l^2 \mu^{(1)}(0)
\end{equation}
and equating those in front of  $\varepsilon^5 {\bf E}$ first yields
\begin{eqnarray}\label{eq. kdv in BR}
-2icK\partial_T \widetilde{A} & = & -\frac{1}{24} \partial_l^4 \mu^{(1)}(0) K^4 \widetilde{A} \\&& - \frac{1}{2} \partial_l^2 \mu^{(1)}(0) K^2 \int_{-\frac{1}{2\varepsilon}}^\frac{1}{2\varepsilon} \widetilde{A}(T,K-M)\widetilde{A}(T,M)\; dM.
\nonumber
\end{eqnarray}
%resp.
%\begin{equation}\label{eq. KdV2*}
%\partial_T \widetilde{A} = -\frac{1}{48c} \partial_l^4 w^{(1)}(0) iK^3 \widetilde{A} - \frac{1}{4c} \partial_l^2 w^{(1)}(0) iK \int_{-\frac{1}{2\varepsilon}}^\frac{1}{2\varepsilon} \widetilde{A}(T,K-M)\widetilde{A}(T,M)\; dM. 
%\end{equation} 
If we now let the parameter $\varepsilon\to 0$  and $\widetilde{A} \to \hat{A}$ as $\varepsilon\to 0$ , we  formally  obtain the KdV equation for $\hat{A} $ in  Fourier space
\begin{equation}\label{eq. KdV in Fourier}
\partial_T \hat{A} (T,K)= -\frac{1}{48c} \partial_l^4 \mu^{(1)}(0) iK^3 \hat{A} (T,K)- \frac{1}{4c} \partial_l^2 \mu^{(1)}(0) iK (\hat{A}\ast\hat{A})(T,K)
\end{equation}
or 
\begin{equation}\label{eq. KdV1}
\partial_T {A}(T,X) = \frac{1}{48c} \partial_l^4 \mu^{(1)}(0) \partial_X^3 {A}(T,X) - \frac{1}{4c} \partial_l^2 \mu^{(1)}(0) \partial_X (A^2)(T,X)
\end{equation}
for $ A $ in physical space.
%which is consistent to the KdV equation \eqref{eq. KdV1} in physical space.
%\begin{remark}{\rm 

There is now a problem. The original system lives in Bloch space but 
the KdV equation lives in Fourier space. Therefore, we set  
\begin{equation} \label{KdVBloch}
\widetilde{A}(T,K) = \chi_{[-\frac{1}{4}\delta_0,\frac{1}{4}\delta_0 ]} (l) \hat{A}(T,K) 
\end{equation}
for the transition from Bloch to Fourier space and vice versa, where we 
extend $\widetilde{A}$ periodically to $\R$ by the above periodicity conditions \eqref{period1}. 
%}\end{remark}

\section{The improved approximation}

\label{secresidual}

\label{Subkap. The improved approach}
The residual 
\[ \widetilde{\text{Res}}_2 (\widetilde{U})(t,l,x)= \partial_t^2 \widetilde{U}(t,l,x) + \widetilde{\mathcal{B}}_l^2\widetilde{U}(t,l,x)+\widetilde{\mathcal{B}}_l^2 N(\widetilde{U})(t,l,x), \]
contains all terms which remain after inserting the above ansatz into 
\eqref{eq. Bouss1*}. 
In order to estimate the error made by the KdV approximation with our subsequent approach
 the residual in Bloch space should not be bigger than  $\mathcal{O}(\varepsilon^7)$. 
 To achieve this, we make an improved ansatz that eliminates all terms up to  order $\mathcal{O}(\varepsilon^6)$. 
To do this, we add higher order terms to the original ansatz, i.e., we make the modified ansatz
\[ \widetilde{U} = 
\varepsilon^2 \widetilde{\Psi}  = 
\widetilde{V}_\text{KdV}(t,l) f^{(1)}(l,x) + \widetilde{U}^\perp_\text{KdV} (t,l,x) \]
with $\widetilde{V}_\text{KdV}$ from \eqref{eq. Ansatz3 V}  and
\eqref{KdVBloch} as well as 
\begin{equation}\label{eq. Ansatz3 U}
\widetilde{U}^\perp_\text{KdV}(t,l,x) = \varepsilon^3 \widetilde{B}_1(T,K,x) \,{\bf E}(t,l) + \varepsilon^5 \widetilde{B}_2 (T,K,x) \,{\bf E}(t,l)
\end{equation}
and insert this into the second equation of \eqref{eq. Bouss2*}. 
With some abuse of notation, we write for instance 
$ K \widetilde{B}_1 = K \widetilde{B}_1(T,K,x) $. 
We find
\begin{align*}
&(\varepsilon^{9} \partial_T^2 \widetilde{B}_1 + \varepsilon^{11} \partial_T^2 \widetilde{B}_2)\, {\bf E} - 2icK (\varepsilon^7 \partial_T \widetilde{B}_1 + \varepsilon^9 \partial_T \widetilde{B}_2)\, {\bf E} - c^2 K^2 ( \varepsilon^5 \widetilde{B}_1 + \varepsilon^7 \widetilde{B}_2)\,{\bf E} \\
&= -( \varepsilon^3 \widetilde{\mathcal{B}}_l^2 \widetilde{B}_1 + \varepsilon^5 \widetilde{\mathcal{B}}_l^2 \widetilde{B}_2)\,{\bf E} - \left( \varepsilon^2 \widetilde{\mathcal{B}}_l^2((\widetilde{A}f^{(1)})\ast(\widetilde{A}f^{(1)})) \right. \\
&\qquad +\left. 2\varepsilon^4 \widetilde{\mathcal{B}}_l^2 ((\widetilde{A}f^{(1)})\ast\widetilde{B}_1) + 2\varepsilon^6 \widetilde{\mathcal{B}}_l^2((\widetilde{A}f^{(1)})\ast\widetilde{B}_2) \right. \\
&\qquad + \left. \varepsilon^6 \widetilde{\mathcal{B}}_l^2 (\widetilde{B}_1\ast\widetilde{B}_1) + 2\varepsilon^8 \widetilde{\mathcal{B}}_l^2 (\widetilde{B}_1\ast\widetilde{B}_2) + \varepsilon^{10} \widetilde{\mathcal{B}}_l^2 (\widetilde{B}_2\ast\widetilde{B}_2) \right) {\bf E}\\
&\qquad + \varepsilon^4 \frac{1}{2}\partial_l^2 \omega^{(1)}(0) K^2 
 \chi_{[-\frac{1}{2}\delta_0,\frac{1}{2}\delta_0]}(\varepsilon K) \langle f^{(1)},(\widetilde{A}f^{(1)})\ast(\widetilde{A}f^{(1)})\rangle_{L_\Gamma^2}f^{(1)} \, {\bf E} + \mathcal{O}(\varepsilon^7),
\end{align*}
where all convolution terms gain one extra $\varepsilon$-power due to the substitution $m=\varepsilon M$
in the convolution integrals such as  the one  above \eqref{eq. beta}. We now choose $\widetilde{B}_1$ and $\widetilde{B}_2$ such that the $\mathcal{O}(\varepsilon^3)$- resp. $\mathcal{O}(\varepsilon^5)$-terms cancel, i.e., we set
\begin{equation}\label{eq. B1}
\widetilde{B}_1 (T,K,x) = -((\widetilde{A}f^{(1)})\ast(\widetilde{A}f^{(1)}))(T,K,x)
\end{equation}
and
\begin{align*}
c^2 K^2 \widetilde{B}_1 = \,\,&\widetilde{\mathcal{B}}_l^2 \widetilde{B}_2 + 2\widetilde{\mathcal{B}}_l^2 ((\widetilde{A}f^{(1)})\ast\widetilde{B}_1) 
\\
& - \frac{1}{2}\partial_l^2 \omega^{(1)}(0) K^2  \chi_{[-\frac{1}{2}\delta_0,\frac{1}{2}\delta_0]}(\varepsilon K) \langle f^{(1)},(\widetilde{A}f^{(1)})\ast(\widetilde{A}f^{(1)})\rangle_{L_\Gamma^2}f^{(1)}. 
\end{align*}
Because of \eqref{eq. c}  and \eqref{eq. B1} we can choose $\widetilde{B}_2$ as follows
\begin{equation}\label{eq. B2}
\begin{array}{rcl}
\widetilde{B}_2(T,K,x) &=& 2 ((\widetilde{A}f^{(1)})\ast(\widetilde{A}f^{(1)})\ast(\widetilde{A}f^{(1)}))(T,K,x)\\[1mm] & &
 - \displaystyle \frac{1}{2}\partial_l^2 \omega^{(1)}(0) K^2 \widetilde{\mathcal{B}}_l^{-2} P_l^{s}((\widetilde{A}f^{(1)})\ast(\widetilde{A}f^{(1)}))(T,K,x),
\end{array}
\end{equation}

where we use that $\widetilde{\mathcal{B}}_l^2$ is invertible on the range of $P_l^{s}$.

\section{Analytic properties of the approximation}

In this section we collect a number of analytic properties  of the  approximation  \eqref{eq. Ansatz3 U}.

%\begin{remark}{\rm
Due to the $L^2$-scaling property
$ \| A(\varepsilon^{-1}\cdot ) \|_{L_2} = \varepsilon^{1/2} \| A( \cdot  ) \|_{L_2} $
we lose a factor $\varepsilon^{1/2}$ in  $\widetilde{\mathcal{H}}^2$. 
Because of the $L^1$-scaling property and the scaling property of Fourier transform 
this is not the case  in 
the space 
\[ \widetilde{\mathcal{C}}^2 := \{ \widetilde{U}\in L^1(\T_1,L_\Gamma^2): \; \widetilde{u}_j \in L^1(\T_1,H^2(I_{0,j})),\; j\in\{ 0,+,- \},\; \eqref{eq. KH1*}-\eqref{eq. KH2*}\;\text{satisfied} \} \]
equipped with the norm
\[ \| \widetilde{U} \|_{\widetilde{\mathcal{C}}^2} = \int_{-\frac{1}{2}}^\frac{1}{2} \left( \| \widetilde{u}_0(l,\cdot) \|_{H^2(I_{0,0})} + \| \widetilde{u}_+(l,\cdot) \|_{H^2(I_{0,+})} + \| \widetilde{u}_-(l,\cdot) \|_{H^2(I_{0,-})} \right)\; dl ,\]
which we will use in the following to estimate the KdV approximation in the equations for the error.
%}\end{remark}

For estimating the KdV approximation and later the residual  
we follow \cite[Subsection 5.3]{GPS16} and  introduce a weight with respect to the variable $l$, namely
\[ \rho_{\varepsilon,s} (l) =\left( 1+ \left(\frac{l}{\varepsilon}\right)^2 \right)^\frac{s}{2},\quad s\geq 0. \]
In the following we use that  Fourier transform is an isomorphism between the Sobolev space $H^s$ and $L^2$ with weight $\rho_{1,s}$.

The subsequent assumption $A\in C([0,T_0],H^6)$ 
for the solutions of the KdV equation implies that  $\hat{A}\rho_{1,6} \in L^2$ and thus $\hat{A}\rho_{1,5} \in L^1$ due to
\[ \| \hat{A}\rho_{1,5} \|_{L^1} \leq \| \hat{A}\rho_{1,6} \|_{L^2} \| \rho_{1,-1} \|_{L^2} \leq C \| \hat{A}\rho_{1,6} \|_{L^2}. \]
Moreover, since  $f^{(1)}_j$ is smooth with respect to $x$ for $j\in\{ 0,+,- \}$, it is also bounded in $H^2(I_{0,j})$ due to the compactness of  $I_{0,j}$.

This implies that the main part 
$\widetilde{V}_\text{kdv} f^{(1)}  $
of the approximation $ \varepsilon^2 \widetilde{\Psi}  $ is 
$ \mathcal{O}(\varepsilon^2) $-bounded in the space $   \widetilde{\mathcal{C}}^2 $.
Using this estimate and the choices $\widetilde{B}_1$ and $\widetilde{B}_2$ from \eqref{eq. B1} and \eqref{eq. B2} gives that $\widetilde{U}^\perp_\text{kdv}  $  is of order $\mathcal{O}(\varepsilon^4)$ in $ \widetilde{\mathcal{C}}^2$.
Therefore,  we obtain that the approximation $ \varepsilon^2 \widetilde{\Psi}$  satisfies
\begin{equation}\label{eq. Psibound}
\sup_{t\in[0,T_0/ \varepsilon^3]} \| \varepsilon^2 \widetilde{\Psi} \|_{\widetilde{\mathcal{C}}^2} \leq C_{\Psi} \varepsilon^{2}
\end{equation}
with an $\varepsilon$-independent constant $C_{\Psi} >0$. 

%Subsequently we have to estimate products of the approximation and the 
%error for which we use estimates, cf. 
%\cite[Subsection 5.3]{GPS16}, like 
%\[ \| (\widetilde{V}\ast\widetilde{W})\rho_{\varepsilon,s} \|_{\widetilde{\mathcal{H}}^2} \leq C\| \widetilde{V}\rho_{\varepsilon, s} \|_{\widetilde{\mathcal{C}}^2} \| \widetilde{W}\rho_{\varepsilon,s} \|_{\widetilde{\mathcal{H}}^2}, \quad s\geq 0, \]
%which is based on 
%Young's inequality. 

\section{Estimates for the residual}

With  our  choice  of the approximation  and since the equations are even w.r.t. the Bloch wave number $ l $, all terms up to  order $\mathcal{O}(\varepsilon^6)$ vanish. Thus the residual is of 
formal order $\mathcal{O}(\varepsilon^7)$ in Bloch space and of formal order $\mathcal{O}(\varepsilon^8)$ in physical space.  In Sobolev norms we obtain  
\begin{lemma}\label{Lem. Resabsch}
Let $A\in C([0,T_0],H^6)$ be a solution of the KdV equation \eqref{eq. KdV1} for a $T_0>0$. Then there exist an $\varepsilon_0>0$ and a constant $C_\text{Res}>0$ such that in Bloch space  
\begin{equation}\label{eq. Resabsch1}
\begin{split}
\sup_{t\in [0,T_0 / \varepsilon^3]} \| \widetilde{\text{Res}}_2 (\varepsilon^2 \widetilde{\Psi}) \|_{\widetilde{\mathcal{H}}^2} \leq C_\text{Res} \varepsilon^{15/2} 
\end{split}
\end{equation}
and in physical space 
\begin{equation}\label{eq. Resabsch2}
\begin{split}
\sup_{t\in [0,T_0 / \varepsilon^3]} \| \text{Res}_2 (\varepsilon^2 \Psi) \|_{\mathcal{H}^2} \leq C_\text{Res} \varepsilon^{15/2}
\end{split}
\end{equation}
for all $\varepsilon\in (0,\varepsilon_0)$.
\end{lemma}
\noindent {\bf Proof.} 
The rigorous proof of the smallness of the residual goes along the lines of  the proof of \cite[Lemma 5.5]{GPS16} 
and uses estimates such as the subsequent Lemma \ref{Lem. Help}
to turn the formal smallness into smallness in function spaces.
Since this is straightforward and documented in a various papers, we 
refrain from giving the details. In this way the first estimate is proven. 
The second estimate then follows from Lemma \ref{Lem. Isom}.
\qed

Since $ \mathcal{H}^2 \subset \mathcal{H}^1$ as Banach spaces, 
the second estimate holds in $ \mathcal{H}^1$, too.

\begin{remark}{\rm
a) Due to the $L^2$-scaling property
$\| A(\varepsilon^{-1}\; \cdot \;) \|_{L_2} = \varepsilon^{1/2} \| A(\; \cdot \; ) \|_{L_2}$,
we lose a factor $\varepsilon^{1/2}$ if we estimate the residual in the space $\widetilde{\mathcal{H}}^2$. 

b)  In the proof of Lemma \ref{Lem. Resabsch}, one has to estimate the difference between the curve of eigenvalues
$ \mu^{(1)}(l)$ and its fourth Taylor polynomial at $l=0$. Since $ \mu^{(1)}$ 
is even, we need for this estimate six derivatives of $A$ and the following lemma \cite[Lemma 5.4]{GPS16} for $m=0$ and $s=6$: 
\begin{lemma}\label{Lem. Help}
Let $m,s\geq 0$ and $g:\T_1\rightarrow \R$ be a function with
\[ |g(l)|\leq C|l|^s,\quad l\in\T_1. \]
Then, we have
\[ \| \rho_{1,m}(\cdot) g(\cdot)\widetilde{A}(\varepsilon^{-1}(\cdot)) \|_{L^2(\T_1)} \leq C \varepsilon^{s+1/2} \| \rho_{1,m+s} \hat{A} \|_{L^2(\R)}. \]
\end{lemma}
}
\end{remark}

Lemma \ref{Lem. Resabsch} also gives  an estimate for the difference between $\varepsilon^2 \Psi$ from Section \ref{secresidual} and 
\[\varepsilon^2 \Psi_\text{KdV}:=\varepsilon^2 A( \varepsilon (x-c t) ,\varepsilon^3 t) f_1(0,x).\] 
Applying  Lemma \ref{Lem. Resabsch}  to the expansion of $ f_1 $ w.r.t. $ l $ which can be found in Lemma \ref{Lem. Eigfkt} and to the appearing cut-off functions gives
\begin{lemma}
Let $A\in C([0,T_0],H^6)$ be a solution of the KdV equation \eqref{eq. KdV1} for $T_0>0$. Then there exists  an $\varepsilon$-independent  constant $C>0$, such that
\begin{equation}\label{eq. Ungl2}
\sup_{t\in[0,T_0/ \varepsilon^3]} \| \varepsilon^2 \Psi - \varepsilon^2 \Psi_{\text{KdV}} \|_{L^\infty} \leq C \varepsilon^{5/2}.
\end{equation}
\end{lemma}

Subsequently, we will also need 
\begin{lemma}\label{Lem. Resabsch77}
Let $A\in C([0,T_0],H^6)$ be a solution of the KdV equation \eqref{eq. KdV1} for a $T_0>0$. Then there exist an $\varepsilon_0>0$ and a constant $C_\text{Res}>0$ such that 
%in Bloch space  
for all $\varepsilon\in (0,\varepsilon_0)$:
%\begin{equation}\label{eq. Resabsch1}
%\begin{split}
%\sup_{t\in [0,T_0 / \varepsilon^3]} \| \mathcal{A}^{-1} \widetilde{\text{Res}}_2 (\varepsilon^2 \widetilde{\Psi}) \|_{\widetilde{\mathcal{H}}^2} \leq C_\text{Res} \varepsilon^{13/2},
%\end{split}
%\end{equation}
%and in physical space that
\begin{equation}\label{eq. Resabsch3}
\begin{split}
\sup_{t\in [0,T_0 / \varepsilon^3]} \| \mathcal{B}^{-1} \text{Res}_2 (\varepsilon^2 \Psi) \|_{\mathcal{H}^2} \leq C_\text{Res} \varepsilon^{13/2}.
\end{split}
\end{equation}
\end{lemma}
It is obvious that we can apply $ \mathcal{B}^{-1}  $ on those residual terms
coming from the right-hand side of the 
Boussinesq equation since all these terms have a $ \mathcal{B}^2 $ in front.
Hence, it remains to control the terms coming from the left-hand side of the 
Boussinesq equation. The operator $ \mathcal{B}^{-1}  $ has a singularity 
$ \mathcal{O}(1/l) $ at the wave number $ l = 0 $ and at $ \omega = 0 $.
Therefore, we have to guarantee that in Bloch space 
the terms coming from $ \partial^2_T A $ possess a representation 
which shows that at the wave number $ l = 0 $ and at $ \omega = 0 $ they  
are at least proportional to $ \mathcal{O}(l) $.
Expressing the time derivatives of $ A $ through the right-hand side 
of the KdV equation gives such a representation.
In fact, following  the arguments of the proof of \cite[Lemma 2.1]{BDS19} yields  that $ \partial^2_T A $ can be expressed 
as spatial derivative, i.e. in Fourier or Bloch space that 
at the wave number $ l = 0 $ and at $ \omega = 0 $ this term  is at least proportional to $ \mathcal{O}(l) $.

\section{Error estimates}

\label{sec8}

Our approximation result is as follows
\begin{theorem}\label{Main theorem}
Let $A\in C([0,T_0],H^8)$ be a solution of the KdV equation \eqref{eq. KdV1} for $T_0>0$. Then there exist $\varepsilon_0>0$ and $C>0$ such that for all $\varepsilon\in (0,\varepsilon_0)$ there exists a solution $U\in C([0,T_0/\varepsilon^3],\mathcal{H}^1)$ of the Boussinesq equation \eqref{eq. Bouss1} with
\begin{equation}\label{eq. Main theorem}
\sup_{t \in [0,T_0/\varepsilon^3]} \sup_{x \in \R} |U(x,t)-\varepsilon^2 A( \varepsilon (x-c t) ,\varepsilon^3 t)f^{(1)}(0,x) | \leq C \varepsilon^{5/2}.
\end{equation}
\end{theorem}
\noindent
{\bf Proof.}
We write  a solution $U$ of  \eqref{eq. Bouss2} as  sum of  the approximation $\varepsilon^2 \Psi$ and an the error  $\varepsilon^{7/2} R$, i.e.,
\begin{equation}\label{eq. Ansatz Fehler}
U=\varepsilon^2 \Psi+\varepsilon^{7/2} R.
\end{equation}
Inserting \eqref{eq. Ansatz Fehler} into \eqref{eq. Bouss2} gives an equation for the error
\begin{equation}\label{eq. Fehlergleichung}
\partial_t^2 R = -\mathcal{B}^2(\partial_x) R - 2\varepsilon^2 \mathcal{B}^2(\partial_x) (\Psi R) - \varepsilon^{7/2} \mathcal{B}^2 (\partial_x) (R^2) + \varepsilon^{-7/2} \text{Res}_2 (\varepsilon^2 \Psi),
\end{equation}
where as before $R=(R_0,R_+,R_-)$ and $\Psi=(\psi_0,\psi_+,\psi_-)$ are multiplied
component-wise. 
The only term which makes difficulties to bound the error function $ R $ on the long 
$ \mathcal{O}(1/\varepsilon^3)$-time scale is the term $ - 2\varepsilon^2 \mathcal{B}^2(\Psi R) $.
In order to control this term, we construct 
an energy
$\mathcal{E}$ with $\frac{d}{dt} \mathcal{E}=\mathcal{O}(\varepsilon^3)$. 
Before we do so, we make the following remark.
\begin{remark} \label{Bem. Resabsch}{\rm 
In the following energy estimates the operator $ (I + \mathcal{A}^2) $ will fall on the residual terms. Therefore, we define
$$
\textrm{Res}_1 = (I + \mathcal{A}^2)\textrm{Res}_2.
$$
Obviously, we can use the estimates from Lemma \ref{Lem. Resabsch} and 
Lemma \ref{Lem. Resabsch77} for $\textrm{Res}_2$ also for $\textrm{Res}_1$ by increasing the regularity by two, i.e., we have to assume $ A \in  C([0, T_0], H^8) $.
}\end{remark}
In order to find this energy we first take the $\mathcal{L}^2$-scalar product of 
$\partial_t(I + \mathcal{A}^2)R$ with \eqref{eq. Fehlergleichung}. 
We obtain
\begin{equation} \label{eq. Fehlergleichung*}
\left\langle \partial_t^2 R,\partial_t(I + \mathcal{A}^2)R \right\rangle_{\mathcal{L}^2} 
= - s_1 - 2\varepsilon^2 s_2 - \varepsilon^{7/2}s_3 +\varepsilon^{-7/2} s_4,
\end{equation}
where
\begin{eqnarray*}
s_1 &= & \left\langle \mathcal{B}^2 R,\partial_t(I + \mathcal{A}^2)R \right\rangle_{\mathcal{L}^2} ,\\
s_2 &= & \left\langle \mathcal{B}^2(\Psi R),\partial_t(I + \mathcal{A}^2)R \right\rangle_{\mathcal{L}^2} ,\\
s_3 &= & \left\langle \mathcal{B}^2(R^2),\partial_t(I + \mathcal{A}^2)R \right\rangle_{\mathcal{L}^2} ,\\
s_4 &= & \left\langle \text{Res}_2 (\varepsilon^2 \Psi),\partial_t(I + \mathcal{A}^2)R \right\rangle_{\mathcal{L}^2}.
\end{eqnarray*}
Using again $ \|\partial_x R \|_{\mathcal{L}^2} = \|\mathcal{A} R \|_{\mathcal{L}^2} $ the term on the left hand side can be written as 
\begin{align*}
\left\langle \partial_t^2 R,\partial_t(I+ \mathcal{A}^2)R \right\rangle_{\mathcal{L}^2} = \frac{1}{2} \frac{d}{dt} \left\| \partial_t R \right\|_{\mathcal{L}^2}^2 + \frac{1}{2} \frac{d}{dt} \left\| \partial_t \partial_x R \right\|_{\mathcal{L}^2}^2.
\end{align*}
For the terms on the right hand side we find
\begin{eqnarray*}
s_1 &= & \left\langle  \mathcal{A}^2 R,\partial_t R \right\rangle_{\mathcal{L}^2} = \left\langle  \mathcal{A} R,\partial_t  \mathcal{A} R \right\rangle_{\mathcal{L}^2} = \frac{1}{2} \frac{d}{dt} \left\|  \partial_x R \right\|_{\mathcal{L}^2}^2, \\
s_2 &= & \left\langle \mathcal{A}^2 (\Psi R),\partial_t R \right\rangle_{\mathcal{L}^2} = \left\langle \partial_x (\Psi R),\partial_t \partial_x R \right\rangle_{\mathcal{L}^2} \\
&=& \frac{1}{2} \frac{d}{dt} \left\langle \partial_x R,\Psi(\partial_x R) \right\rangle_{\mathcal{L}^2} - \varepsilon \left\langle \partial_x R,(\partial_\tau \Psi)(\partial_x R) \right\rangle_{\mathcal{L}^2} \\
&&\qquad - \varepsilon \left\langle \partial_t R,(\partial_X\Psi) (\partial_x R) \right\rangle_{\mathcal{L}^2} - \varepsilon^2 \left\langle \partial_t R,(\partial_X^2 \Psi) R \right\rangle_{\mathcal{L}^2}, \\
s_3 &= &\left\langle \mathcal{A}^2 (R^2),\partial_t R \right\rangle_{\mathcal{L}^2} = \left\langle \partial_x (R^2),\partial_t \partial_x R \right\rangle_{\mathcal{L}^2} \\
&=& \frac{d}{dt}\left\langle R,(\partial_x R)^2 \right\rangle_{\mathcal{L}^2} - \left\langle \partial_x R,(\partial_t R) (\partial_x R) \right\rangle_{\mathcal{L}^2}, \\
s_4 & = & \left\langle (1-\partial_x^2)\text{Res}_2 (\varepsilon^2 \Psi),\partial_t R \right\rangle_{\mathcal{L}^2} 
 = \left\langle \text{Res}_1 (\varepsilon^2 \Psi),\partial_t R \right\rangle_{\mathcal{L}^2},
\end{eqnarray*}
where $\tau=\varepsilon t$ and where all
boundary terms at the vertices cancel due to the 
conditions \eqref{eq. KH1*} and \eqref{eq. KH2*}.

Next we take the $\mathcal{L}^2$-scalar product of 
$\partial_t\mathcal{B}^{-2}R$ with \eqref{eq. Fehlergleichung}. We find 
\begin{equation}\label{eq. Fehlergleichung2*}
\left\langle \partial_t^2 R,\partial_t\mathcal{B}^{-2}R \right\rangle_{\mathcal{L}^2} 
= - s_5 - 2\varepsilon^2  s_6- \varepsilon^{7/2} s_7 + \varepsilon^{-7/2} s_8,
\end{equation}
where
\begin{eqnarray*}
s_5 &= & \left\langle \mathcal{B}^2 R,\partial_t\mathcal{B}^{-2}R \right\rangle_{\mathcal{L}^2} ,\\
s_6 &= &\left\langle \mathcal{B}^2(\Psi R),\partial_t\mathcal{B}^{-2}R \right\rangle_{\mathcal{L}^2} ,\\
s_7 &= & \left\langle \mathcal{B}^2(R^2),\partial_t\mathcal{B}^{-2}R \right\rangle_{\mathcal{L}^2} ,\\
s_8 &= &  \left\langle \text{Res}_2 (\varepsilon^2 \Psi),\partial_t\mathcal{B}^{-2}R \right\rangle_{\mathcal{L}^2}.
\end{eqnarray*}
The term on the left hand side can be written as 
%\begin{align*}
%\left\langle \partial_t^2 R,\partial_t\mathcal{B}^{-2}R \right\rangle_{\mathcal{L}^2} = 
%\left\langle \partial_t^2 \mathcal{B}^{-1} R,\partial_t\mathcal{B}^{-1}R \right\rangle_{\mathcal{L}^2} = \frac{1}{2} \frac{d}{dt} \left\| \partial_t  \mathcal{B}^{-1} R \right\|_{\mathcal{L}^2}^2 
%\end{align*}
\begin{align*}
\left\langle \partial_t^2 R,\partial_t\mathcal{B}^{-2}R \right\rangle_{\mathcal{L}^2} =  \frac{1}{2} \frac{d}{dt} \left\| \mathcal{B}^{-1} \partial_t  R \right\|_{\mathcal{L}^2}^2 
%+ \frac{1}{2} \frac{d}{dt} \left\| \partial_t \partial_x^{-1}R \right\|_{\mathcal{L}^2}^2 
\end{align*}
%since formally $\mathcal{B}^{-2}=I+\mathcal{A}^{-2}$. 
For the terms on the right hand side we find
\begin{eqnarray*}
s_5 &= & \left\langle R,\partial_t R \right\rangle_{\mathcal{L}^2} = \frac{1}{2}\frac{d}{dt} \| R \|_{\mathcal{L}^2}^2, \\
s_6 &= & \left\langle \Psi R,\partial_t R \right\rangle_{\mathcal{L}^2} = \frac{1}{2}\frac{d}{dt}\left\langle R,\Psi R \right\rangle_{\mathcal{L}^2} - \varepsilon \left\langle R,( \partial_\tau \Psi ) R \right\rangle_{\mathcal{L}^2}, \\
s_7 &= & \left\langle R^2,\partial_t R \right\rangle_{\mathcal{L}^2} = \frac{1}{3}\frac{d}{dt}\left\langle R , R^2 \right\rangle_{\mathcal{L}^2}, \\
s_8 &= & \left\langle \text{Res}_2 (\varepsilon^2 \Psi),\partial_t \mathcal{B}^{-2} R \right\rangle_{\mathcal{L}^2} = \left\langle \mathcal{B}^{-1} \text{Res}_2 (\varepsilon^2 \Psi),\partial_t \mathcal{B}^{-1} R \right\rangle_{\mathcal{L}^2}.
\end{eqnarray*}
We collect all terms which have a time derivative in front  in the energy $\mathcal{E}$, i.e., we set
\begin{align*}
\mathcal{E} &= \frac{1}{2} \left\| \partial_t R \right\|_{\mathcal{L}^2}^2 + \frac{1}{2} \left\| \partial_t \partial_x R \right\|_{\mathcal{L}^2}^2 + \frac{1}{2} \left\| \partial_x R \right\|_{\mathcal{L}^2}^2 + \frac{1}{2} \| \mathcal{B}^{-1} \partial_t R \|_{\mathcal{L}^2}^2 
%+ \frac{1}{2} \| \partial_t R \|_{\mathcal{L}^2}^2 
+ \frac{1}{2} \| R \|_{\mathcal{L}^2}^2 \\
&\qquad+ \varepsilon^2 \left\langle \partial_x R,\Psi(\partial_x R) \right\rangle_{\mathcal{L}^2} + \varepsilon^{7/2}\left\langle R,(\partial_x R)^2 \right\rangle_{\mathcal{L}^2} + \varepsilon^2 \left\langle R,\Psi R \right\rangle_{\mathcal{L}^2} + \frac{1}{3}\varepsilon^{7/2} \left\langle R , R^2 \right\rangle_{\mathcal{L}^2}.
\end{align*}
For sufficiently small $\varepsilon>0$, we can estimate the second line by the first one. Thus for all $M>0$ there exist constants $C_1,\varepsilon_1 > 0$, such that for all $\varepsilon\in (0,\varepsilon_1)$ we have
\[ \| R \|_{\mathcal{H}^1} \leq C_1 \mathcal{E}^{1/2} \]
as long as $\mathcal{E}\leq M$ holds.
We estimate the remaining terms on the right hand side in terms of $\mathcal{E}$ by using the Cauchy-Schwarz inequality and Sobolev's embedding theorem.
We find
\begin{align*}
\left|\left\langle \partial_x R,(\partial_\tau \Psi)(\partial_x R) \right\rangle_{\mathcal{L}^2}\right| &\leq \| \partial_\tau \Psi \|_{L^\infty} \| \partial_x R \|_{\mathcal{L}^2}^2 \leq C\mathcal{E}, \\
\left|\left\langle \partial_t R,(\partial_X\Psi) (\partial_x R) \right\rangle_{\mathcal{L}^2}\right| &\leq \| \partial_X \Psi \|_{L^\infty} \| \partial_t R \|_{\mathcal{L}^2} \| \partial_x R \|_{\mathcal{L}^2} \leq C\mathcal{E}, \\
\left|\left\langle \partial_t R,(\partial_X^2 \Psi) R \right\rangle_{\mathcal{L}^2}\right| &\leq \| \partial_X^2 \Psi \|_{\mathcal{L}^2} \| \partial_t R \|_{\mathcal{L}^2} \| R \|_{L^\infty} \leq C  \varepsilon^{-1/2}\mathcal{E}, \\
\left|\left\langle \partial_x R,(\partial_t R) (\partial_x R) \right\rangle_{\mathcal{L}^2}\right| &\leq \| \partial_t R \|_{L^\infty} \| \partial_x R \|_{\mathcal{L}^2}^2 
\\&\leq C\sqrt{\| \partial_t R \|_{\mathcal{L}^2}\| \partial_t \partial_x R \|_{\mathcal{L}^2}} \| \partial_x R \|_{\mathcal{L}^2}^2 \leq C\mathcal{E}^{3/2}, \\
\left|\left\langle \text{Res}_1 (\varepsilon^2 \Psi),\partial_t R \right\rangle_{\mathcal{L}^2}\right| &\leq \| \text{Res}_1 (\varepsilon^2 \Psi) \|_{\mathcal{L}^2} \|\partial_t R\|_{\mathcal{L}^2}\\ &\leq C\varepsilon^{15/2} \|\partial_t R\|_{\mathcal{L}^2} \leq C \varepsilon^{15/2} \mathcal{E}^{1/2}, \\
\left| \left\langle R,( \partial_\tau \Psi ) R \right\rangle_{\mathcal{L}^2} \right| &\leq \| \partial_\tau \Psi \|_{L^\infty} \| R \|_{\mathcal{L}^2}^2 \leq C\mathcal{E}, \\
\left| \left\langle \mathcal{B}^{-1} \text{Res}_2 (\varepsilon^2 \Psi),\partial_t \mathcal{B}^{-1} R \right\rangle_{\mathcal{L}^2} \right| &\leq \| \mathcal{B}^{-1} \text{Res}_2 (\varepsilon^2 \Psi) \|_{\mathcal{L}^2}\| \partial_t \mathcal{B}^{-1} R \|_{\mathcal{L}^2} \leq C \varepsilon^{13/2} \mathcal{E}^{1/2},
\end{align*}
where we used $ \| \partial_\tau \Psi \|_{L^\infty} + \| \partial_X \Psi \|_{L^\infty} + \varepsilon^{1/2} \| \partial_X^2 \Psi \|_{\mathcal{L}^2} \leq C \|  \Psi \|_{\mathcal{C}^2} $ and 
$ \| R \|_{L^\infty}  \leq C \|R \|_{\mathcal{H}^1} $,  in particular in the first three lines.

By using $\mathcal{E}^{1/2}\leq 1+\mathcal{E}$ we find that 
the energy $\mathcal{E}$ satisfies the inequality
%and $ \| \partial_\tau \Psi \|_{L^\infty}, \| \partial_X \Psi \|_{L^\infty}, \| \partial_X^2 \Psi \|_{L^\infty} < \infty$
\begin{align*}
\frac{d}{dt} \mathcal{E} \leq C \varepsilon^3 \mathcal{E} + C\varepsilon^{7/2} \mathcal{E}^{3/2} + C \varepsilon^3.
\end{align*}
Under the assumption that  $C\varepsilon^{1/2}\mathcal{E}^{1/2} \leq 1$, Gronwall's inequality finally yields
\[ \sup_{t\in [0,T_0/\varepsilon^3]} \mathcal{E}(t) = CT_0 e^{(C+1)T_0}=:M = \mathcal{O}(1) \]
and thus 
\[ \sup_{t\in [0,T_0/\varepsilon^3]} \| R(t)\|_{\mathcal{H}^1} = \mathcal{O}(1). \]
If we choose $\varepsilon_2>0$ so small that $C\varepsilon_2^{1/2}M^{1/2} \leq 1$ holds,  
we obtain the required estimate for all $\varepsilon\in (0,\varepsilon_0)$ with $\varepsilon_0=\min(\varepsilon_1,\varepsilon_2)$.  \qed
\label{secenergy}

\section{Discussion}

\label{sec9}

We close this paper with two remarks. 
\begin{remark}{\rm
There are no serious obstacles to do the same analysis 
for the Boussinesq model posed on other periodic quantum graphs
as long as the spectrum is of the qualitative form displayed in Figure 
\ref{Fig. SpektrumPlot} near $ (l,\omega) = (0,0) $.
The theory for instance also applies to the Boussinesq equation posed on the graphs 
displayed in \cite[Section 7]{GPS16}, in particular it applies 
to the necklace graph where the upper and lower circles have different length.
Also for PDEs on the real line there is no abstract theorem which applies 
to all dispersive systems. 
%Even there the validity results are example oriented.
The variety of problems on quantum graphs is much bigger since now the equation and the quantum graph can vary. Therefore, a general KdV validity theory
on quantum graphs would be extremely abstract and of minor use. 
}\end{remark}
\begin{remark}{\rm As in \cite{GPS16} at points  $ (l,\omega) \neq  (0,0) $,
without at the intersection points of the spectral curves at $ l = 0 $ and $ \omega \neq 0 $, 
the NLS equation can be derived through a multiple scaling ansatz.
At the intersection points a Dirac equation can be derived.
However, the justification of the NLS approximation
for the Boussinesq equation posed on the necklace graph is so far out of reach
due to the quadratic terms.
In the justification analysis these have to be eliminated 
through normal form transformations. 
However, due to the conditions at the vertex points, in the symmetric case, where
$ u_+ = u_- $ and where the necklace graph can be identified with the real line, 
the derivative makes a jump at the periodically distributed vertex points. Therefore, we are in a very 
non-smooth situation and the existing theory for the validity of the NLS 
approximation on the real line for problems with spatially periodic coefficients 
and quadratic terms does not apply, cf. \cite{BSTU06,BCCLS08}.
Moreover, an additional difficulty occurs due to the zero eigenvalues at 
$ k = 0 $ which leads to some resonance. See \cite{DS06} 
for the handling of this resonance and other possible resonances 
for some Boussinesq equation on the real line.
In \cite{GPS16} the NLS justification  has only  been handled for cubic nonlinearities.  
There, no normal form transformations are necessary.
}\end{remark}

\subsection*{Conflict of interest statement}
The authors declare that they do not have any conflict of interest.

\bibliographystyle{alpha}

\bibliography{literatur}

\newcommand{\etalchar}[1]{$^{#1}$}
\begin{thebibliography}{GMWZ14}

\bibitem[BBC{\etalchar{+}}08]{BCCLS08}
Carsten Blank, Martina~Chirilus Bruckner, Christopher Chong, Vincent Lescarret,
  Guido Schneider, and Hannes Uecker.
\newblock A remark about the justification of the nonlinear {Schr{\"o}dinger}
  equation in quadratic spatially periodic media.
\newblock {\em Z. Angew. Math. Phys.}, 59(3):554--557, 2008.

\bibitem[BDS19]{BDS19}
Roman {Bauer}, Wolf-Patrick {D\"ull}, and Guido {Schneider}.
\newblock {The Korteweg-de Vries, Burgers and Whitham limits for a spatially
  periodic Boussinesq model}.
\newblock {\em {Proc. R. Soc. Edinb., Sect. A, Math.}}, 149(1):191--217, 2019.

\bibitem[BK13]{Kuchment}
Gregory Berkolaiko and Peter Kuchment.
\newblock {\em Introduction to quantum graphs}, volume 186 of {\em Mathematical
  Surveys and Monographs}.
\newblock American Mathematical Society, Providence, RI, 2013.

\bibitem[BSTU06]{BSTU06}
Kurt Busch, Guido Schneider, Lasha Tkeshelashvili, and Hannes Uecker.
\newblock Justification of the nonlinear {Schr{\"o}dinger} equation in
  spatially periodic media.
\newblock {\em Z. Angew. Math. Phys.}, 57(6):905--939, 2006.

\bibitem[CCPS12]{CCPS12}
Martina {Chirilus-Bruckner}, Christopher {Chong}, Oskar {Prill}, and Guido
  {Schneider}.
\newblock {Rigorous description of macroscopic wave packets in infinite
  periodic chains of coupled oscillators by modulation equations}.
\newblock {\em {Discrete Contin. Dyn. Syst., Ser. S}}, 5(5):879--901, 2012.

\bibitem[{Cra}85]{Cr85}
Walter {Craig}.
\newblock {An existence theory for water waves and the Boussinesq and Korteweg-
  deVries scaling limits}.
\newblock {\em {Commun. Partial Differ. Equations}}, 10:787--1003, 1985.

\bibitem[DJ89]{Drazin}
Philip~G. Drazin and Robin~S. Johnson.
\newblock {\em Solitons: an introduction}.
\newblock Camb. Texts Appl. Math. Cambridge etc.: Cambridge University Press,
  1989.

\bibitem[DS06]{DS06}
Wolf-Patrick D{\"u}ll and Guido Schneider.
\newblock Justification of the nonlinear {Schr{\"o}dinger} equation for a
  resonant {Boussinesq} model.
\newblock {\em Indiana Univ. Math. J.}, 55(6):1813--1834, 2006.

\bibitem[D{\"u}l12]{Du12}
Wolf-Patrick D{\"u}ll.
\newblock {Validity of the Korteweg-de Vries approximation for the
  two-dimensional water wave problem in the arc length formulation}.
\newblock {\em {Commun. Pure Appl. Math.}}, 65(3):381--429, 2012.

\bibitem[GMWZ14]{GMWZ14}
Jeremy {Gaison}, Shari {Moskow}, J.~Douglas {Wright}, and Qimin {Zhang}.
\newblock {Approximation of polyatomic FPU lattices by KdV equations}.
\newblock {\em {Multiscale Model. Simul.}}, 12(3):953--995, 2014.

\bibitem[GPS16]{GPS16}
Steffen Gilg, Dmitry Pelinovsky, and Guido Schneider.
\newblock Validity of the {NLS} approximation for periodic quantum graphs.
\newblock {\em NoDEA Nonlinear Differential Equations Appl.}, 23(6):Art. 63,
  30, 2016.

\bibitem[GSU04]{Gallay}
Thierry Gallay, Guido Schneider, and Hannes Uecker.
\newblock Stable transport of information near essentially unstable localized
  structures.
\newblock {\em Discrete Contin. Dyn. Syst. Ser. B}, 4(2):349--390, 2004.

\bibitem[GSU21]{GSU20}
Steffen Gilg, Guido Schneider, and Hannes Uecker.
\newblock Modulated waves on graphene like metric graphs.
\newblock {\em Mathematische Nachrichten}, accepted, 2021.

\bibitem[KN86]{KN86}
Tadayoshi {Kano} and Takaaki {Nishida}.
\newblock {A mathematical justification for Korteweg-de Vries equation and
  Boussinesq equation of water surface waves}.
\newblock {\em {Osaka J. Math.}}, 23:389--413, 1986.

\bibitem[Pel11]{Pelinovsky}
Dmitry~E. Pelinovsky.
\newblock {\em Localization in periodic potentials}, volume 390 of {\em London
  Mathematical Society Lecture Note Series}.
\newblock Cambridge University Press, Cambridge, 2011.
\newblock From Schr\"{o}dinger operators to the Gross-Pitaevskii equation.

\bibitem[PS17]{Pelinovsky*}
Dmitry Pelinovsky and Guido Schneider.
\newblock Bifurcations of standing localized waves on periodic graphs.
\newblock {\em Ann. Henri Poincar\'{e}}, 18(4):1185--1211, 2017.

\bibitem[RS80]{Reed}
Michael Reed and Barry Simon.
\newblock {\em Methods of modern mathematical physics. {I}}.
\newblock Academic Press, Inc. [Harcourt Brace Jovanovich, Publishers], New
  York, second edition, 1980.
\newblock Functional analysis.

\bibitem[SU17]{SU17book}
Guido Schneider and Hannes Uecker.
\newblock {\em Nonlinear PDEs: a dynamical systems approach}.
\newblock American Mathematical Society, 2017.

\bibitem[SW00a]{SW00equa}
Guido {Schneider} and C.~Eugene {Wayne}.
\newblock {Counter-propagating waves on fluid surfaces and the continuum limit
  of the Fermi-Pasta-Ulam model}.
\newblock In {\em {International conference on differential equations.
  Proceedings of the conference, Equadiff '99, Berlin, Germany, August 1--7,
  1999. Vol. 1}}, pages 390--404. Singapore: World Scientific, 2000.

\bibitem[SW00b]{SW00}
Guido {Schneider} and C.~Eugene {Wayne}.
\newblock {The long-wave limit for the water wave problem. I: The case of zero
  surface tension}.
\newblock {\em {Commun. Pure Appl. Math.}}, 53(12):1475--1535, 2000.

\bibitem[SW02]{SW02}
Guido {Schneider} and C.~Eugene {Wayne}.
\newblock {The rigorous approximation of long-wavelength capillary-gravity
  waves}.
\newblock {\em {Arch. Ration. Mech. Anal.}}, 162(3):247--285, 2002.

\end{thebibliography}

\end{document}